# An algorithm based on DQM with modified trigonometric cubic B-splines for solving coupled viscous Burgers' equations

Brajesh Kumar Singh*, Pramod Kumar

Department of Applied Mathematics, School for Physical Sciences, Babasaheb Bhimrao Ambedkar University Lucknow - 226025 INDIA
bksingh0584@gmail.com

**Abstract:** This paper deals with a new algorithm called modified trigonometric cubic B-spline differential quadrature method (MTB-DQM) for numerical computation of the time dependent partial differential equations. Specially, the numerical computation of the Burgers' equation is obtained by using MTB-DQM with time integration algorithm. The MTB-DQM is new DQM with modified trigonometric cubic B-splines as basis. The initial boundary value system of Burgers' equation is first transformed into an initial value system of ordinary differential equations (ODEs) by means of MTB-DQM. The resulting system of ODEs is solved by using an optimal five stage four order strong stability preserving Runge–Kutta method (SSP-RK54). The accuracy and efficiency of the method is illustrated by five test problems in terms of $L_2$ and $L_\infty$ error norms and their comparisons with existing results. The MTB-DQM produces better results than the results due to almost all the existing schemes. The MTB-DQM is shown conditionally stable using the matrix stability analysis method for various grid points.

## 1. Introduction

The Burgers' equation, Navier–Stoke's equation without the stress term, is the easiest nonlinear physical model various flows problems consisting of hydrodynamic turbulence, sound and shock wave theory, vorticity transportation, wave processes in thermo-elastic medium, dispersion in porous media, mathematical modeling of turbulent fluid, continuous stochastic processes. This equation was first introduced by Bateman [1], for details study, we refer readers to [2-6].

Consider initial valued system of coupled viscous Burgers' equation in dimension (1+1) as:

$$\begin{cases} \dfrac{\partial u}{\partial t} = \dfrac{\partial^2 u}{\partial x^2} - \eta u \dfrac{\partial u}{\partial x} - \alpha \dfrac{\partial (uv)}{\partial x}, \\ \dfrac{\partial u}{\partial t} = \dfrac{\partial^2 v}{\partial x^2} - \xi v \dfrac{\partial v}{\partial x} - \beta \dfrac{\partial (uv)}{\partial x}, \\ u(x,0) = \phi(x), \quad v(x,0) = \psi(x), \quad x \in \Omega_1, \ \mathrm{t} > 0. \end{cases} \tag{1.1}$$

with the Dirichlet boundary conditions

$$u(a,t) = g_1(t), \ \ u(b,t) = g_2(t), \ \ v(a,t) = g_3(t), v(b,t) = g_4(t), \ \ t > 0, \tag{1.2}$$

and (2+1) dimensional nonlinear coupled viscous Burger's equation:

$$\begin{cases} \dfrac{\partial u}{\partial t} + u \dfrac{\partial u}{\partial x} + v \dfrac{\partial u}{\partial y} = \upsilon \left( \dfrac{\partial^2 u}{\partial x^2} + \dfrac{\partial^2 u}{\partial y^2} \right), \\ \dfrac{\partial v}{\partial t} + u \dfrac{\partial v}{\partial x} + v \dfrac{\partial v}{\partial y} = \upsilon \left( \dfrac{\partial^2 v}{\partial x^2} + \dfrac{\partial^2 v}{\partial y^2} \right), \\ u(x,y,0) = \phi(x,y) \text{ and } v(x,y,0) = \psi(x,y), \quad (x,y) \in \Omega, \ t > 0 \end{cases} \tag{1.3}$$

with the Dirichlet boundary conditions as follows:

$$u(x,y,t) = \xi(x,y,t) \ \text{ and } \ v(x,y,t) = \zeta(x,y,t), \quad (x,y) \in \partial\Omega, \ t > 0, \tag{1.4}$$

Where $\Omega_1 = [a, b]$, $\Omega = \{(x,y) : a \le x \le b, \ c \le x \le d\} \rightarrow$ computational domains in one and two dimension, and $\partial\Omega$ is the boundary of $\Omega$ , $\xi$, $\eta$ real constants and $\alpha, \beta$ are arbitrary constants depending on the system parameters such as Peclet number, stokes velocity of particles due to gravity and Brownian diffusivity [7]. $u(x,y,t)$, $v(x,y,t) \rightarrow$ velocity components, $\dfrac{\partial u}{\partial t} \rightarrow$ unsteady term, $u\dfrac{\partial u}{\partial x} \rightarrow$ nonlinear convection term, $\upsilon\left(\dfrac{\partial^2 u}{\partial x^2} + \dfrac{\partial^2 u}{\partial y^2}\right) \rightarrow$ diffusion term and $\upsilon \rightarrow$ coefficient of viscosity ($\upsilon > 0$) and $\phi$, $\psi$, $g_1, g_2, g_3, g_4$, $\xi$ and $\zeta$ are known functions.

In the recent years, a lot of efforts have been made for the computation of accuracy and efficiency of numerical schemes for Burgers' equation with various values of kinematic viscosity. Burgers' equation has been solved using various analytical and numerical schemes such as Hofe Cole transformation [5-6], finite element method [8], finite difference method [9,32], implicit finite difference method [10,30], Quartic B-spline collocation method [31], compact finite difference method [11-13], Fourier Pseudospectral method [14], variational iteration method [15], cubic B-spline collocation scheme [16-17], modified cubic B-splines collocation method [18], modified trigonometric cubic B-spline collocation method [19], extended B-spline collocation method [20], reproducing kernel function method [21], quadratic B-spline finite elements [22]. For more schemes on two dimensional Buregers' equation, the interested readers are referred to [57-62]. One dimensional coupled viscous Burgers'equation has been solved numerically using various techniques, among others, lattice Boltzmann method [23], differential quadrature method [24], A fully implicit finite difference [25], a composite scheme based on finite difference and Haar wavelets [26], discrete Adomian decomposition method [27], ], variational iteration method [28], algorithms based on cubic spline function technique [29].

The differential quadrature method (DQM), developed by Bellman et al. [33], is widely used for numerical solution of partial differential equations (PDEs). After the seminal paper of Bellman et al. [33], and Quan and Chang [34-35] DQM has been implemented with various type of set of basis functions, among others, cubic B-spline differential quadrature methods [36], DQM based on Fourier expansion and Harmonic function [37-38], sinc DQM [39], generalized DQM [40] and modified cubic B-spline differential quadrature method [41-44], modified extended cubic B-spline differential quadrature (mECDQ) method [45], Polynomial based DQM [53], quartic B-spline based DQM [52,54], Quartic and quintic B-spline methods [55], exponential cubic B-spline DQM [56].

The B-splines (piece-wise smooth polynomials) are capable to handle local phenomena, and so, they have more influence in comparison to other set of basis functions. The main goal of this paper is to develop a new method referred to as "modified trigonometric cubic-B-spline differential quadrature method (MTB-DQM)" for numerical computation of the time dependent partial differential equations. Specially, the numerical computation of the Burgers' equation is done using MTB-DQM with time integration algorithm. MTB-DQM is a DQM with modified trigonometric cubic-B-splines as a new set of basis functions. The initial boundary valued system of coupled viscous Burgers' equation is first converted into an initial value system of ordinary differential equations (ODEs) by means of MTB-DQM. The SSP-RK54 scheme is used to solve the resulting system of ODEs. The SSP-RK54 scheme is chosen due to it's reduce storage space which results in less accumulation of the numerical errors. The accuracy and efficiency of the method is illustrated by five test problems in terms of $L_2$ and $L_\infty$ error norms and their comparisons with existing results. The MTB-DQM produces better results than the results due to almost all the existing schemes. The MTB-DQM is shown conditionally stable using the matrix stability analysis method for various grid points. The MTB-DQM solutions of Burgers' equation are computed without transforming the equation and without any linearization technique.

## 2. Description: modified trigonometric B-spline differential quadrature method

This section deals with the description of MTB-DQM for Burgers' equation. The weighting coefficients being depends on only the grid spacing, the domains $\Omega$ and $\Omega_1$ defined by

$\Omega_1 = \{x \in \mathbb{R} : a \le x \le b\}$, and $\Omega = \{(x, y) \in \mathbb{R}^2 : a \le x \le b,\ c \le y \le d\}$

are portioned uniformly in each direction with the following knots:

$a = x_1 < x_2, \cdots < x_i < \cdots < x_{N_x - 1} < x_{N_x} = b,$

$c = y_1 < y_2, \cdots < y_j < \cdots < y_{N_y - 1} < y_{N_y} = d,$

Where $h_x = \dfrac{b-a}{N_x - 1}$; $h_y = \dfrac{d-c}{N_y - 1}$ denote the space step in $x$, $y$ directions, respectively. Let

$x_i \in \Omega_1, (x_i, y_j) \in \Omega$ be the generic grid points and

$u_i := u_i(t) = u(x_i, t);$

$u_{ij} := u_{ij}(t) = u(x_i, y_j, t).$

The trigonometric cubic B-spline function [19] is defined as follows

$$\mathrm{T}_m(x) = \frac{1}{\omega}\begin{cases} p^3(x_m), & x \in [x_m, x_{m+1}) \\ p(x_m)(p(x_m)q(x_{m+2}) + q(x_{m+3})p(x_{m+1})) + q(x_{m+4})p^2(x_{m+1}), & x \in [x_{m+1}, x_{m+2}) \\ q(x_{m+4})(p(x_{m+1})q(x_{m+3}) + q(x_{m+4})p(x_{m+2})) + p(x_m)q^2(x_{m+3}), & x \in [x_{m+2}, x_{m+3}) \\ q^4(x_{m+4}), & x \in [x_{m+3}, x_{m+4}) \end{cases} \tag{2.1}$$

where $p(x_m) = \sin\left(\dfrac{x - x_m}{2}\right)$, $q(x_m) = \sin\left(\dfrac{x_m - x}{2}\right)$, $\omega = \sin\left(\dfrac{h_x}{2}\right)\sin(h_x)\sin\left(\dfrac{3h_x}{2}\right)$.

where the set $\{T_0, T_1, \ldots, T_{N_x}, T_{N_x+1}\}$ of trigonometric cubic B-spline functions form a basis for the computational region $\Omega_1$. Setting

$$a_1 = \frac{\sin^2\left(\dfrac{h_x}{2}\right)}{\sin(h_x)\sin\left(\dfrac{3h_x}{2}\right)}, \quad a_2 = \frac{2}{1 + 2\cos(h_x)}, \quad a_4 = \frac{3}{4\sin\left(\dfrac{3h_x}{2}\right)} = -a_3,$$

$$\mathrm{a}_5 = \frac{3 + 9\cos(h_x)}{16\sin^2\left(\dfrac{h_x}{2}\right)\left(2\cos\left(\dfrac{h_x}{2}\right) + \cos\left(\dfrac{3h_x}{2}\right)\right)}, \quad a_6 = -\frac{3\cos^2\left(\dfrac{h_x}{2}\right)}{\sin^2\left(\dfrac{h_x}{2}\right)(2 + 4\cos(h_x))}$$

The functional value of trigonometric cubic B-spline $T_{ij} := T_i\left(x_j\right)$ and its first and second derivatives $T'_{ij} := T'_i\left(x_j\right),\ T''_{ij} := T''_i\left(x_j\right)$ at $x_j$ is, respectively, read as:

$$T_{ij} = \begin{cases} a_2, & i-j=0 \\ a_1, & i-j=\pm 1 \\ 0, & \text{otherwise} \end{cases};\qquad T'_{ij} = \begin{cases} a_4, & i-j=1 \\ a_3, & i-j=-1 \\ 0, & \text{otherwise} \end{cases};\qquad T''_{ij} = \begin{cases} a_6, & i-j=0 \\ a_5, & i-j=\pm 1 \\ 0, & \text{otherwise} \end{cases} \tag{2.2}$$

The modified trigonometric cubic B-splines are defined in similar manner as in [41, 45]:

$$\left.\begin{array}{l} \sigma_1(x) = T_1(x) + 2T_0(x), \\ \sigma_2(x) = T_2(x) - T_0(x) \\ \sigma_m(x) = T_m(x), \qquad\qquad m = 3,\ldots,N_x - 2 \\ \sigma_{N-1}(x) = T_{N-1}(x) - T_{N+1}(x), \\ \sigma_N(x) = T_N(x) + 2T_{N+1}(x) \end{array}\right\}, \tag{2.3}$$

where $\left\{\sigma_1,\ \sigma_2,\ \ldots,\sigma_{N_x}\right\}$ forms a basis of $\Omega_1$.

The r[th] order partial derivatives $(r \geq 2)$ of $u\left(x,t\right)$ with respect to $x$ approximated at $x_i$ read:

$$\frac{\partial^r u_i}{\partial x^r} = \sum_{k=1}^{N} a_{ik}^{(r)} u_k,\ i \in \Delta_{N_x} \tag{2.4}$$

The r[th] order partial derivatives $(r \geq 2)$ of $u(x,y,t)$ with respect to $x,\ y$ at the grid point $(x_i, y_j)$ can be computed as follows:

$$\begin{cases} \dfrac{\partial^r u}{\partial x^r}\left(x_i,y_j\right) = \displaystyle\sum_{k=1}^{N} a_{ik}^{(r)} u_{kj}, & i \in \Delta_{N_x} \\ \dfrac{\partial^r u}{\partial y^r}\left(x_i,y_j\right) = \displaystyle\sum_{k=1}^{M} b_{jk}^{(r)} u_{ik}, & j \in \Delta_{N_y} \end{cases},\qquad \begin{cases} \dfrac{\partial^r v}{\partial x^r}\left(x_i,y_j\right) = \displaystyle\sum_{k=1}^{N} a_{ik}^{(r)} v_{kj}, & i \in \Delta_{N_x} \\ \dfrac{\partial^r v}{\partial y^r}\left(x_i,y_j\right) = \displaystyle\sum_{k=1}^{M} b_{jk}^{(r)} v_{ik}, & j \in \Delta_{N_y} \end{cases} \tag{2.5}$$

where $a_{ij}^{(r)}$ and $b_{ij}^{(r)}$ $\left(r = 1, 2\right)$ are the weighting coefficients of the $r^{th}$ order partial derivatives with respect to $x$ and $y$.

**2.1 Computation of the weighting coefficients**

In order to evaluate the weighting coefficients $a_{ij}^{(1)}$ of Eq. (2.4), the modified extended cubic B-spline $\sigma_m\left(x\right), m \in \Delta_{N_x}$ are used. Setting $\sigma'_{mi} := \phi'_m\left(x_i\right), \sigma_{mi} := \sigma_m\left(x_i\right)$. Accordingly, the approximation for the first order spatial derivative is given by:

$$\sigma_{mi} = \sum_{\ell=1}^{N_x} a_{i\ell}^{(1)} \sigma_{m\ell}, \quad i,m \in \Delta_{N_x}. \tag{2.6}$$

After setting, $\Sigma = [\sigma_{m\ell}]$, $A = [a_{i\ell}^{(1)}]$, and $\Sigma' = [\sigma'_{mi}]$, Eq. (2.6) can be written as the set of tridiagonal system of linear equations as follows:

$$\Sigma A^T = \Sigma'. \tag{2.7}$$

The coefficient matrix $\Sigma$ of order $N_x$ can be computed from Eq. (2.2) and Eq.(2.3) as:

$$\Sigma = \begin{bmatrix} (a_2 + 2a_1) & a_1 & & & & & \\ 0 & a_2 & a_1 & & & & \\ & a_1 & a_2 & a_1 & & & \\ & & \ddots & \ddots & \ddots & & \\ & & & a_1 & a_2 & a_1 & \\ & & & & a_1 & a_2 & 0 \\ & & & & & a_1 & (a_2 + 2a_1) \end{bmatrix}$$

and in particular the column of the matrix $\Sigma'$ read:

$$\Sigma'[1] = \begin{bmatrix} a_4 \\ a_3 - a_4 \\ 0 \\ 0 \\ \vdots \\ 0 \\ 0 \end{bmatrix}, \Sigma'[2] = \begin{bmatrix} a_4 \\ 0 \\ a_3 \\ 0 \\ \vdots \\ 0 \\ 0 \end{bmatrix}, \Sigma'[3] = \begin{bmatrix} 0 \\ a_4 \\ 0 \\ a_3 \\ 0 \\ \vdots \\ 0 \end{bmatrix}, \cdots, \Sigma'[N_x - 1] = \begin{bmatrix} 0 \\ 0 \\ \vdots \\ 0 \\ a_4 \\ 0 \\ a_3 \end{bmatrix}, \Sigma'[N_x] = \begin{bmatrix} 0 \\ 0 \\ \vdots \\ 0 \\ 0 \\ a_4 \\ a_3 - a_4 \end{bmatrix}.$$

Noticed that the tridiagonal system (2.7) has unique solution because $\Sigma$ is invertible. *Thomas algorithm* [47] is used to solve the system (2.7), which yields the weighting coefficients $a_{i\ell}^{(1)}, i, \ell \in \Delta_{N_x}$. Similarly, the weighting coefficients $b_{i\ell}^{(1)}$ can be computed by considering the grid descritization in $y$ direction.

. The existence of various basis spanning of the problem domain, the weighting coefficients ($a_{i\ell}^{(2)}, i, \ell \in \Delta_{N_x}$) of the second-order derivatives can be determined using various set of basis functions. One can determine these weighting coefficients using the second order spatial derivative approximation as [36]:

$$\sigma''_{mi} = \sum_{\ell=1}^{N} a_{i\ell}^{(2)} \sigma_{m\ell}, \quad i, m \in \Delta_{N_x}. \tag{2.8}$$

Similar to system (2.7), the above system (2.8) yields $a_{i\ell}^{(2)}, i, \ell \in \Delta_{N_x}$. The existence of more than one basis functions to span an $N_x$-dimensional vector space gives an opportunity to compute the weighting coefficients in the same space with other set of basis functions. We prefer Shu's *r*th order ($r \geq 2$) recursive formula [40, 42] based of polynomial based DQM, to compute $a_{i\ell}^{(2)}, i, \ell \in \Delta_{N_x}$; $b_{i\ell}^{(2)}, i, \ell \in \Delta_{N_y}$ as follows:

$$\begin{cases} a_{i\ell}^{(r)} = r\left( a_{i\ell}^{(1)} a_{ii}^{(r-1)} - \dfrac{a_{i\ell}^{(r-1)}}{x_i - x_\ell} \right), & i \neq \ell, \ i, \ell \in \Delta_{N_x}; \\ a_{ii}^{(r)} = -\displaystyle\sum_{\ell=1, \ell \neq i}^{N} a_{i\ell}^{(r)}, & i = \ell, \quad i, \ell \in \Delta_{N_x}. \\ b_{i\ell}^{(r)} = r\left( b_{i\ell}^{(1)} b_{ii}^{(r-1)} - \dfrac{b_{i\ell}^{(r-1)}}{y_i - y_\ell} \right), & i \neq \ell, \ i, \ell \in \Delta_{N_y}; \\ b_{ii}^{(r)} = -\displaystyle\sum_{\ell=1, \ell \neq i}^{N} b_{i\ell}^{(r)}, & i = \ell, \quad i, \ell \in \Delta_{N_y}. \end{cases}$$

## 2.2 Implementation of the method to the Burger's equation

### a. 1D coupled viscous Burgers' equation

On putting the values of the spatial derivatives approximate using MTB-DQM Eq. (1.1) can be re-written as:

$$\begin{cases} \dfrac{\partial u_i}{\partial t} = \displaystyle\sum_{j=1}^{N_x} a_{ij}^{(2)} u_j - \eta u_i \sum_{j=1}^{N_x} a_{ij}^{(1)} u_j - \alpha \left( u_i \sum_{j=1}^{N_x} a_{ij}^{(1)} v_j + v_i \sum_{j=1}^{N_x} a_{ij}^{(1)} u_j \right), \\ \dfrac{\partial v_i}{\partial t} = \displaystyle\sum_{j=1}^{N_x} a_{ij}^{(2)} v_j - \xi v_i \sum_{j=1}^{N_x} a_{ij}^{(1)} v_j - \beta \left( u_i \sum_{j=1}^{N_x} a_{ij}^{(1)} v_j + v_i \sum_{j=1}^{N_x} a_{ij}^{(1)} u_j \right), \\ u_i(t=0) = \phi(x_i), \\ v_i(t=0) = \psi(x_i), \quad i \in \Delta_{N_x}. \end{cases} \tag{3.1}$$

Keeping boundary conditions (1.2) in mind, Eq. (3.1) reduces to a set of first order ODEs:

$$\begin{cases} \dfrac{\partial u_i}{\partial t} = \sum_{j=1}^{N_x-1} a_{ij}^{(2)} u_j - \eta u_i \sum_{j=1}^{N_x-1} a_{ij}^{(1)} u_j - \alpha \left( u_i \sum_{j=1}^{N_x-1} a_{ij}^{(1)} v_j + v_i \sum_{j=1}^{N_x-1} a_{ij}^{(1)} u_j \right) + F_i, \\ \dfrac{\partial v_i}{\partial t} = \sum_{j=1}^{N_x-1} a_{ij}^{(2)} v_j - \xi v_i \sum_{j=1}^{N_x-1} a_{ij}^{(1)} v_j - \beta \left( u_i \sum_{j=1}^{N_x-1} a_{ij}^{(1)} v_j + v_i \sum_{j=1}^{N_x-1} a_{ij}^{(1)} u_j \right) + G_i, \\ u_i(t=0) = \phi(x_i), \\ v_i(t=0) = \psi(x_i), \quad 1 < i < N_x. \end{cases} \tag{3.2}$$

Setting $\eta_i = \eta u_i, \xi_i = \xi v_i,\ \alpha_i = \alpha u_i,\ \alpha'_i = \alpha v_i, \beta_i = \beta u_i,\ \beta'_i = \beta v_i$, then

$$F_i = \left( a_{i1}^{(2)} u_1 + a_{iN_x}^{(2)} u_{N_x} \right) - \eta_i \left( a_{i1}^{(1)} u_1 + a_{iN_x}^{(1)} u_{N_x} \right) - \alpha_i \left( a_{i1}^{(1)} v_1 + a_{iN_x}^{(1)} v_{N_x} \right) - \alpha'_i \left( a_{i1}^{(1)} u_1 + a_{iN_x}^{(1)} u_{N_x} \right),$$
$$G_i = \left( a_{i1}^{(2)} v_1 + a_{iN_x}^{(2)} v_{N_x} \right) - \xi_i \left( a_{i1}^{(1)} v_1 + a_{iN_x}^{(1)} v_{N_x} \right) - \beta_i \left( a_{i1}^{(1)} v_1 + a_{iN_x}^{(1)} v_{N_x} \right) - \beta'_i \left( a_{i1}^{(1)} u_1 + a_{iN_x}^{(1)} u_{N_x} \right).$$

**b. Two dimensional nonlinear coupled Burger's equation**

Similarly, on implementing MTB-DQM to Eq. (1.3)–(1.4) in space, we get

$$\begin{cases} \dfrac{\partial u_{ij}}{\partial t} = \upsilon \left( \sum_{k=2}^{N_x-1} a_{ik}^{(2)} u_{kj} + \sum_{k=2}^{N_y-1} b_{jk}^{(2)} u_{ik} \right) - u_{ij} \sum_{k=2}^{N_x-1} a_{ik}^{(1)} u_{kj} - v_{ij} \sum_{k=2}^{N_y-1} b_{jk}^{(1)} u_{ik} + F_{ij}, \\ \dfrac{\partial v_{ij}}{\partial t} = \upsilon \left( \sum_{k=2}^{N_x-1} a_{ik}^{(2)} v_{kj} + \sum_{k=2}^{N_y-1} b_{jk}^{(2)} v_{ik} \right) - u_{ij} \sum_{k=2}^{N_x-1} a_{ik}^{(1)} v_{kj} - v_{ij} \sum_{k=2}^{N_y-1} b_{jk}^{(1)} v_{ik} + G_{ij}, \\ u_{ij}(t=0) = \xi(x_i, x_j), v_{ij}(t=0) = \zeta(x_i, x_j), 1 < i < N_x, 1 < j < N_y \end{cases} \tag{3.3}$$

where $u_{ij} = \tau_{ij}$ and $v_{ij} = \kappa_{ij}$ and

$$\begin{aligned} F_{ij} &= \upsilon \left( a_{i1}^{(2)} u_{1j} + a_{iN_x}^{(2)} u_{N_x j} + b_{j1}^{(2)} u_{i1} + b_{jN_y}^{(2)} u_{iN_y} \right) - \tau_{ij} \left( a_{i1}^{(1)} u_{1j} + a_{iN_x}^{(1)} u_{N_x j} \right) - \kappa_{ij} \left( b_{j1}^{(1)} u_{i1} + b_{jN_y}^{(1)} u_{iN_y} \right), \\ G_{ij} &= \upsilon \left( a_{i1}^{(2)} v_{1j} + a_{iN_x}^{(2)} v_{N_x j} + b_{j1}^{(2)} v_{i1} + v_{jN_y}^{(2)} u_{iN_y} \right) - \tau_{ij} \left( a_{i1}^{(1)} v_{1j} + a_{iN_x}^{(1)} v_{N_x j} \right) - \kappa_{ij} \left( b_{j1}^{(1)} v_{i1} + b_{jN_y}^{(1)} v_{iN_y} \right). \end{aligned} \tag{3.4}$$

The initial value systems of first order ODEs as given in Eq. (3.2) or Eq. (3.3) can be solved by numerous techniques, among others, we prefer an *optimal five-stage, order four strong stability-preserving time-stepping Runge–Kutta* (SSP-RK54) scheme [48-49] for numerical computation of initial valued system (3.2) and system (3.3), through the following steps:

$$u^{(1)} = u^m + 0.391752226571890 \Delta t L(u^m)$$

$$u^{(2)} = 0.444370493651235 u^m + 0.555629506348765 u^{(1)} + 0.368410593050371 \Delta t L(u^{(1)})$$

$$u^{(3)} = 0.620101851488403 u^m + 0.379898148511597 u^{(2)} + 0.251891774271694 \Delta t L(u^{(2)})$$

$$u^{(4)} = 0.178079954393132u^m + 0.821920045606868u^{(3)} + 0.544974750228521\Delta tL(u^{(3)})$$

$$u^{(m+1)} = 0.517231671970585u^{(2)} + 0.096059710526147u^{(3)} + 0.063692468666290\Delta tL(u^{(3)})$$
$$+ 0.386708617503269u^{(4)} + 0.226007483236906\Delta tL(u^{(4)})$$

The SSP-RK54 scheme allows low storage and large domain of absolute properties [34-35,52].

## 3. Stability analysis

This section deals with the study of the stability analysis of MTB-DQM for Burgers equation. The stability analysis for both dimension is being similar, we concerned with the stability analysis of two dimensional coupled viscous Burgers' equation (3.3), only. Let $A_2 = \left[a_{ij}^{(2)}\right]; A_1 = \left[a_{ij}^{(1)}\right] \& B_2 = \left[b_{ij}^{(2)}\right]; B_1 = \left[b_{ij}^{(1)}\right]$ be the matrices of the weighting coefficients of order$\left(N_x - 2\right)\left(N_y - 2\right)$. Following [50], setting the term $\alpha_i = \alpha u_i$ in the non linear terms of Eq. (3.2) ,and $u_{ij} = \tau_{ij}$ and $v_{ij} = \kappa_{ij}$ in the non linear termsof Eq. (3.3) are assumed to be locally fixed. Eq. (3.3) can be re-written as

$$\frac{dU}{dt} = \mathrm{B}U + H. \tag{4.1}$$

where

- ✓ $U = \left(u, v\right)^T$ is an unknown vector of the functional values at interior grid points:

  $u = \left(u_{22},\ u_{23}, ..., u_{2(M-1)}, u_{32}, u_{33}, ... u_{3(M-1)}, ... u_{(M-1)2}, u_{(M-1)3} ..., u_{(M-1)(M-1)}\right).$

  $v = \left(v_{22},\ \mathrm{v}_{23}, ..., v_{2(M-1)}, v_{32}, v_{33}, ... v_{3(M-1)}, ... v_{(M-1)2}, v_{(M-1)3}, ..., v_{(M-1)(N-1)}\right).$

- ✓ $O$'s are null matrices.
- ✓ $H = \left(F, G\right)^T$, $F = [F_{ij}]$, $G = [G_{ij}]$, $1 < i < N_x; 1 < j < N_y$ as defined in Eq. (3.5)
- ✓ $B = \begin{bmatrix} A & O \\ O & A \end{bmatrix}$, $A = -\tau_{ij}A_1 - \kappa_{ij}B_1 + \upsilon A_2 + \upsilon B_2$, (4.2)

  where $A_r$ and $B_r$ $(r = 1, 2)$ are square block diagonal matrices (each of order $\left(N_x - 2\right)\left(N_y - 2\right)$) of the weighting coefficients $a_{ij}^{(r)}$, $b_{ij}^{(r)}$, respectively as given below:

$$A_r = \begin{bmatrix} a_{22}^{(r)} I & a_{23}^{(r)} I & \cdots & a_{2(N_x-1)}^{(r)} I \\ a_{32}^{(r)} I & a_{33}^{(r)} I & \cdots & a_{3(N_x-1)}^{(r)} I \\ \vdots & \vdots & \ddots & \vdots \\ a_{(N_x-1)2}^{(r)} I & a_{(N_x-1)3}^{(r)} I & \cdots & a_{(N_x-1)(N_x-1)}^{(r)} I \end{bmatrix};$$

$$B_r = \begin{bmatrix} M_r & O & \cdots & O \\ O & M_r & \cdots & O \\ \vdots & \vdots & \ddots & \vdots \\ O & O & \cdots & M_r \end{bmatrix}, \qquad M_r = \begin{bmatrix} b_{22}^{(r)} & b_{23}^{(r)} & \cdots & b_{2(N_y-1)}^{(r)} \\ b_{32}^{(r)} & b_{33}^{(r)} & \cdots & b_{3(N_y-1)}^{(r)} \\ \vdots & \vdots & \ddots & \vdots \\ b_{(N_y-1)2}^{(r)} & b_{(N_y-1)3}^{(r)} & \cdots & b_{(N_y-1)(N_y-1)}^{(r)} \end{bmatrix}$$

✓ $I$ and $O$ are the matrices of order $(N_y - 2)$ and $(N_x - 2)$, respectively.

The stability of system (4.1) depend on the eigenvalues the matrices B [37]. Moreover, if the solution of the system (4.1) is decreasing in absolute value, it needs that all the eigenvalues of B must have negative real part. The stability region is the set $S = \{z \in C, |R(z)| \le 1, z = \lambda \Delta t\}$ where $R(.)$ is the stability function and $\lambda$ is the eigenvalue of B. The stability region of SSP-RK54 scheme is depicted in Fig 1, see [Fig.5, 64].

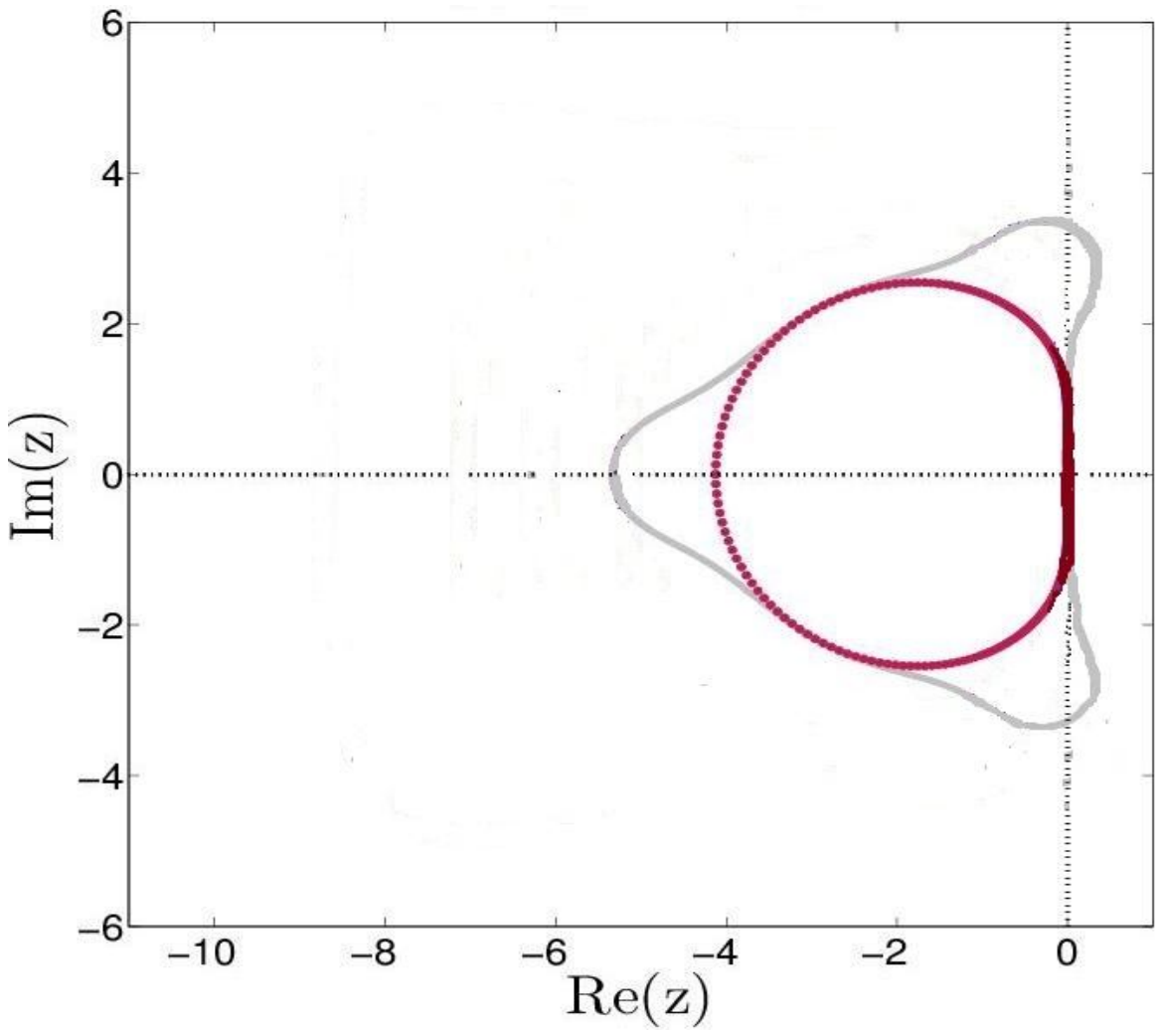

**Fig. 1** Stability region of SSP-RK54 scheme with $z = \lambda \Delta t$

It is sufficient to show the stability of system (4.1) is stable if $\lambda_B \Delta t$ belongs to the stability region of SSP-RK54 scheme for each eigenvalue $\lambda_B$ of B. For more details, see [51, 64].

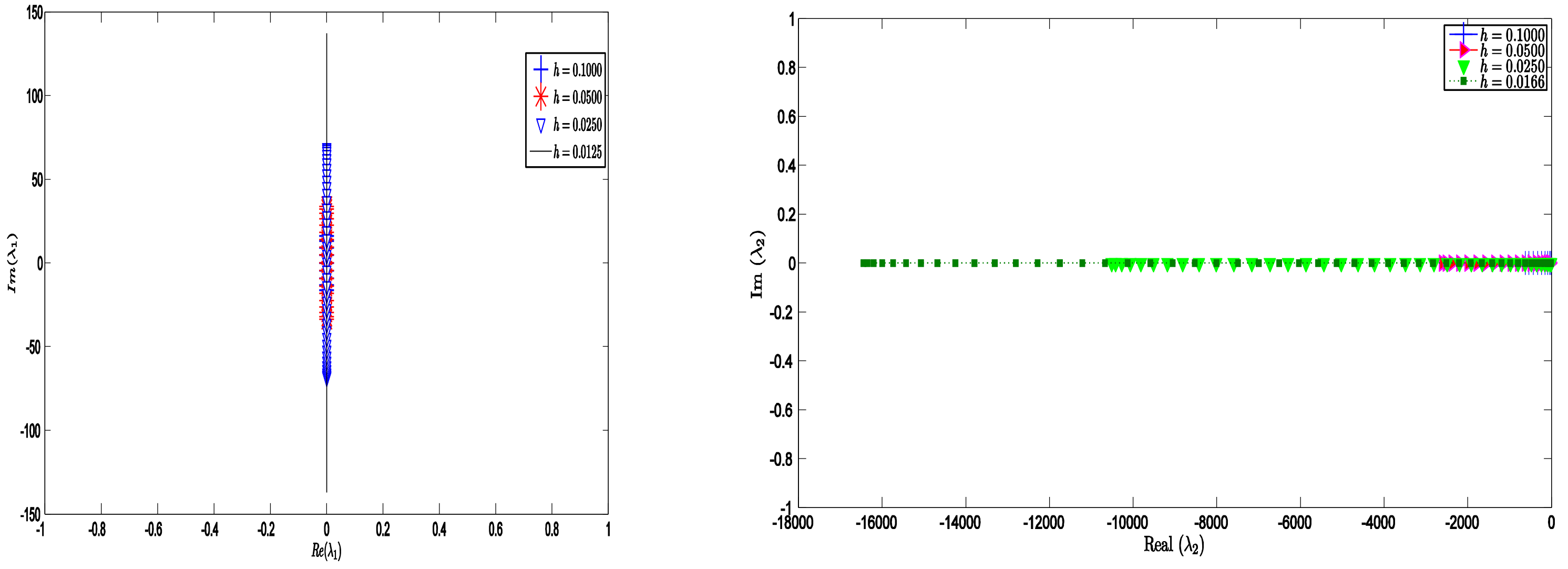

**Fig. 2** Eigen values $\lambda_1$ and $\lambda_2$ for different grid size $h$

It is to be noticed that the eigenvalues of the matrices $A_r$ and $B_r$ $(r = 1, 2)$ have identical nature, and so, it is sufficient to compute the eigenvalues $\lambda_1$ (and $\lambda_2$) of $A_1$ (and $A_2$). The eigenvalues $\lambda_1$ & $\lambda_2$ for different step sizes $h_x = h_y = h$ are depicted in Fig. 2, which confirms that each eigenvalue $\lambda_1$ of the matrices $A_1$ is pure imaginary whereas each eigenvalue $\lambda_2$ of the $A_2$ is real and negative.

It is evident from Fig. 2 and Eq. (4.2) that for a given values of $h$ & $\nu$, there exists $\Delta t$ for which the value $\lambda_B \Delta t$ corresponding to each eignvalue $\lambda_B = 2\nu\lambda_2 - \lambda_1(\tau_0 + \kappa_0)$ of $B$ lies inside the stability region of SSP-RK54. This shows that MTB-DQM produces stable solutions for two dimensional coupled viscous Burger's equation. Similarly, one can demonstrate that MTB-DQM method also produces stable solutions for one dimensional coupled viscous Burger's equation.

## 4. Numerical results and discussion

Now, we consider five test problems of Burgers' equation in one and two dimensions to perform the numerical computation of the proposed MTB-DQM. The accuracy and consistency of the scheme is measured in terms of $L_2$ and $L_\infty$ error norms defined as:

$$L_2 := \sqrt{h\sum_{j=1}^{n}|u_j - u_j^*|^2}; \qquad L_\infty := \max_j |u_j - u_j^*| \tag{5.1}$$

where $u_j$ *and* $u_j^*$ denote exact solution and computed solution at node $x_j$, respectively.

**5.1. 1D coupled Burgers' equation**

**Problem 1** On setting the parameters $\alpha = \beta = 1,\ \xi = \eta = -2$, Eq. (1.1) reduces to

$$\begin{cases} \dfrac{\partial u}{\partial t} = \dfrac{\partial^2 u}{\partial x^2} + 2u\dfrac{\partial u}{\partial x} - \dfrac{\partial(uv)}{\partial x}, & x \in (-\pi,\pi),\ \text{t} > 0, \\ \dfrac{\partial u}{\partial t} = \dfrac{\partial^2 v}{\partial x^2} + 2v\dfrac{\partial v}{\partial x} - \dfrac{\partial(uv)}{\partial x}, & x \in (-\pi,\pi),\ \text{t} \geq 0, \end{cases}$$

The values of $\phi,\ \psi, g_1, g_2, g_3, g_4$ can be obtained from the exact solution as given in [7]:

$u(x,t) = v(x,t) = e^{-t}\sin(x),\ \ x \geq (-\pi,\pi),\ \ t \geq 0.$

The MTB-DQM solution is compared with the results by LBM [23], the order of convergence is reported in Table 1.1. This shows that MTB-DQM has cubic order of convergence in space while LBM [23] is quadratic. Table 1.2 shows that MTB-DQM produces more accurate and smooth results in comparison to FDM and LBM [23]. The comparison of $L_2$ and $L_\infty$ errors in MTB-DQM solutions with the errors from existing schemes at different time levels $t \leq 3,$ with $\Delta t = 0.001$ and $\text{N}_\text{x} = 121$, is reported in Table 1.3.

**Table 1.1**: The comparison of $L_2$ and $L_\infty$ errors and order of convergence $R$ of the Trig-MCB-DQM for $u(x,t)$ Problem 1 with LBM [23] with parameters $\alpha = \beta = 1,\ \xi = \eta = -2$ at $t = 1$

| $N_x$ | MTB-DQM | | | | LBM [53] | | | |
|---|---|---|---|---|---|---|---|---|
| | $L_2$ | $R$ | $L_\infty$ | $R$ | $L_2$ | $R$ | $L_\infty$ | $R$ |
| 10 | 6.69E-03 | | 3.27E-03 | | 3.30E-02 | | 1.15E-02 | |
| 20 | 8.77E-04 | 2.93 | 4.40E-04 | 2.89 | 8.18E-03 | 2.0099 | 3.01E-03 | 1.9375 |
| 40 | 1.03E-04 | 3.08 | 5.13E-05 | 3.10 | 2.01E-03 | 2.0276 | 7.38E-04 | 2.0275 |
| 80 | 8.01E-06 | 3.68 | 4.15E-06 | 3.63 | 4.64E-04 | 2.1123 | 1.71E-04 | 2.1122 |
| 160 | 1.42E-06 | 2.50 | 7.02E-07 | 2.56 | 7.85E-05 | 2.5634 | 2.89E-05 | 2.5638 |

**Table 1.2** The compression of the errors obtained by MCB-DQM, FDM and LBM [23] with the parameters taken as $\Delta t = 0.001$ and $N = 64$

| | $L_2$ | | | $L_\infty$ | | |
|---|---|---|---|---|---|---|
| $t = 1$ | MTB-DQM | LBM [23] | FDM [23] | MTB-DQM | LBM [23] | FDM [23] |
| 0.1 | 8.44E-06 | 3.03E-05 | 8.03E-05 | 2.16E-05 | 2.75E-05 | 7.27E-05 |
| 0.5 | 1.45E-05 | 1.52E-04 | 4.02E-04 | 1.62E-05 | 9.20E-05 | 2.44E-04 |
| 1.0 | 2.02E-05 | 3.03E-04 | 8.03E-04 | 1.02E-05 | 1.12E-04 | 2.95E-04 |
| 2.0 | 3.30E-05 | 6.07E-04 | 1.61E-03 | 4.20E-06 | 8.21E-05 | 2.18E-04 |
| 5.0 | 7.61E-05 | 1.52E-03 | 4.02E-03 | 4.87E-07 | 1.02E-05 | 2.71E-05 |
| 10.0 | 1.50E-04 | 3.04E-03 | 8.06E-03 | 6.64E-09 | 1.38E-07 | 3.66E-07 |
| 20.0 | 3.00E-04 | 6.09E-03 | 1.62E-02 | 6.10E-13 | 1.25E-11 | 3.34E-11 |

**Table 1.3** Comparison with earlier schemes in Problem 1 for $u$ with $\alpha = \beta = 1,\ \xi = \eta = -2\ \ at\ \ t = 1$

| schemes | $N_x$ | $\Delta t$ | $t = 0.5$ $L_\infty$ | $t = 0.5$ $L_2$ | $t = 1$ $L_\infty$ | $t = 1$ $L_2$ | $t = 2$ $L_\infty$ | $t = 3$ $L_\infty$ |
|---|---|---|---|---|---|---|---|---|
| MTB-DQM | 121 | 0.001 | 2.00E-07 | 7.82E-08 | 1.21E-07 | 8.45E-08 | 4.42E-08 | 1.62E-08 |
| DQM[24] | | 0.01 | 1.52E-04 | | 1.84E-04 | | 1.35E-04 | 7.46E-05 |
| CBC[16] | 400 | 0.001 | 6.22E-06 | 1.02E-05 | 7.56E-06 | 2.04E-05 | | |
| CBC[16] | 200 | 0.001 | 4.10E-05 | 2.49E-05 | 8.21E-05 | 0.00003 | $t = 0.1$ | |
| | | | | | | | $L_\infty$ | $L_2$ |
| FPM [14] | 128 | 0.0001 | | | 1.16E-05 | 2.88E-05 | | |
| FFID[25] | 200 | 0.001 | 1.79E-04 | 2.94E-04 | 2.17E-04 | 5.91E-04 | 5.30E-05 | 5.86E-05 |
| CFDH[26] | 16 | 0.001 | 7.27E-06 | 1.05E-08 | 2.38E-05 | 5.03E-07 | 3.77E-05 | 3.27E-08 |
| MTB-DQM | 121 | 0.001 | 2.00E-07 | 7.82E-08 | 1.21E-07 | 8.45E-08 | 3.01E-07 | 7.52E-08 |

It is evident that MTB-DQM produced better solutions as compared to the earlier schemes in [14, 16, 23-26], and found to be in good agreement with the exact solutions. $L_\infty$ and the absolute error norms at different time levels are depicted in Figure 3 whereas the solution behavior is depicted in Figure 4.

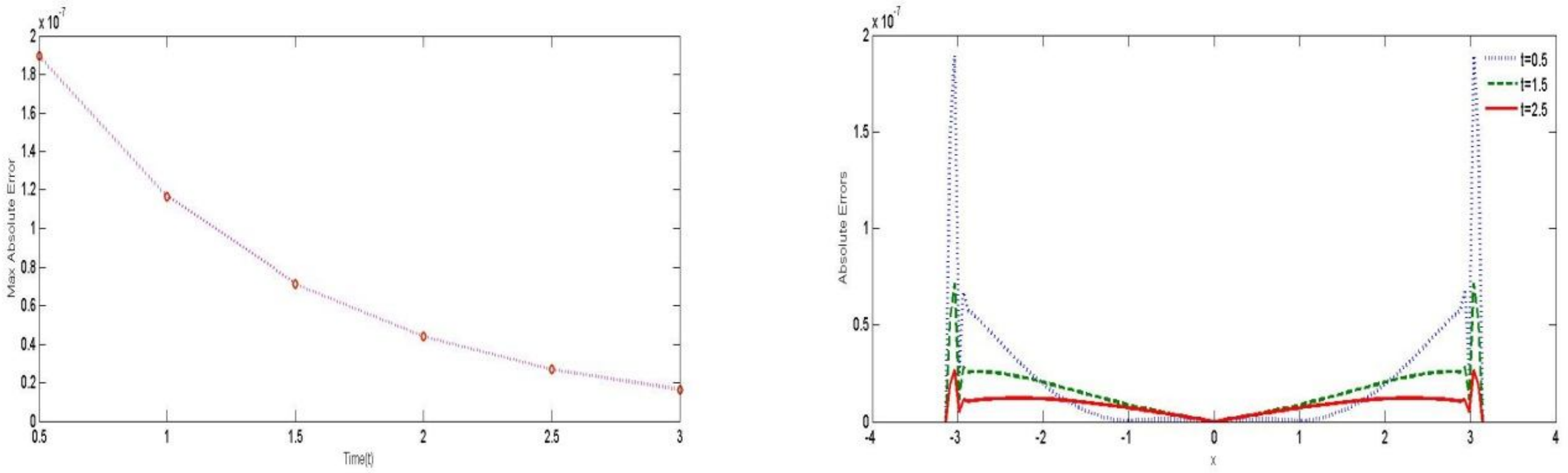

**Figure 3.** The $L_\infty$ and the absolute errors in Example 1 for $u(x,t)$ at different time levels

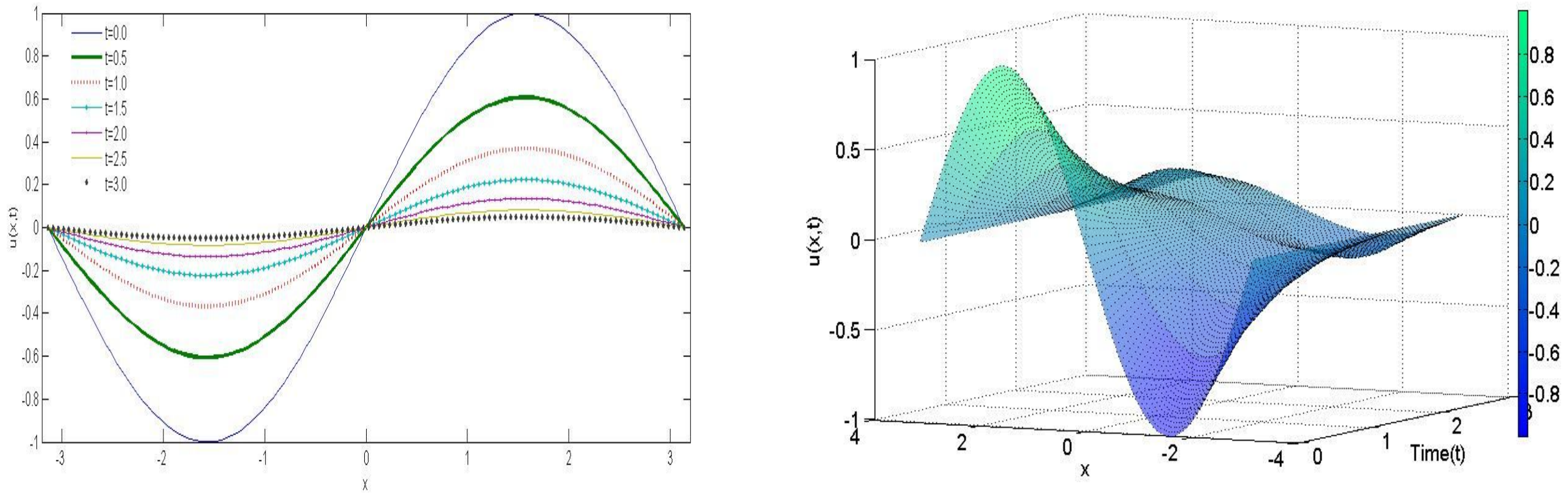

**Figure 4.** The behavior of MTB-DQM solution of Problem 1 for $u$ at $t \le 3.0$ for $\mathrm{N_x} = 121;\ \Delta t = 0.001$

**Problem 2** Consider the 2D Burgers' equations in $\Omega = \left[-0.5, 0.5\right]^2$ as in [27] with $\phi(x,y) = x + y,\ \psi(x,y) = x - y,$ and

$$\left.\begin{array}{ll} u(0,y,t) = \dfrac{y}{1-2t^2}, & u(0.5,y,t) = \dfrac{0.5+y-t}{1-2t^2} \\ v(0,y,t) = \dfrac{-y-2yt}{1-2t^2}, & v(0.5,y,t) = \dfrac{0.5-y-2yt}{1-2t^2} \end{array}\right\},\quad 0 \le y \le 0.5,\quad t \ge 0,$$

$$\left.\begin{array}{ll} u(x,0,t) = \dfrac{x-2xt}{1-2t^2}, & u(x,0.5,t) = \dfrac{x-0.5-2xt}{1-2t^2} \\ v(x,0,t) = \dfrac{x}{1-2t^2}, & v(x,0.5,t) = \dfrac{x-0.5-t}{1-2t^2} \end{array}\right\},\quad 0 \le x \le 0.5,\quad t \ge 0.$$

The exact solutions given in [61] as follows:

$$u(x,y,t) = \frac{x+y-2xt}{1-2t^2},\quad v(x,y,t) = \frac{x-y-2yt}{1-2t^2}.$$

The MTB-DQM solutions are obtained at $t = 0.1$ for uniform mesh grid with $h_x = h_y = 0.025$ and $\Delta t = 10^{-4}$. The MTB-DQM solutions are compared with the exact solutions of Problem 2 with $Re = 1/\upsilon = 80$ for $u$ and $v$ components in Tables 2.1 and 2.2 at various grid points and different time levels $t = 0.1, 0.3, 0.5$. The $L_2,\ L_\infty$ error norms for $u$ and $v$ components are reported in Table 2.3 and Table 2.4 respectively. It is confirmed from the above Tables that the computed results are agreed well with the exact solutions. The physical behavior of MTB-DQM solutions of Problem 2 for $u$ and $v$ components are depicted for $\mathrm{Re}=100$ at $t=0.1$ in Figure 5 whereas for $\mathrm{Re}=80$ at $t=0.5$ in Figure 6.

**Table 2.1** Comparison of the MTB-DQM solution with exact solutions of $u(x,y,t)$ for $\mathrm{Re}=80, \Delta t=10^{-4}$, at different time $t$

| Mesh | $t=0.1$ | | $t=0.3$ | | $t=0.5$ | |
|---|---|---|---|---|---|---|
| | Num. | Exact | Num. | Exact | Num. | Exact |
| 0.1,0.1 | 0.183673 | 0.183673 | 0.170732 | 0.170732 | 0.20004 | 0.20004 |
| 0.3,.01 | 0.346939 | 0.346939 | 0.268259 | 0.268259 | 0.19996 | 0.19996 |
| 0.2,0.2 | 0.367347 | 0.367347 | 0.341465 | 0.341465 | 0.40008 | 0.40008 |
| 0.4,0.2 | 0.530612 | 0.530612 | 0.438991 | 0.438991 | 0.40000 | 0.40000 |
| 0.5,0.3 | 0.714286 | 0.714286 | 0.609723 | 0.609723 | 0.60004 | 0.60004 |
| 0.1,0.3 | 0.387755 | 0.387755 | 0.414670 | 0.414670 | 0.60020 | 0.60020 |
| 0.3,0.4 | 0.653061 | 0.653061 | 0.634166 | 0.634166 | 0.80020 | 0.80020 |
| 0.2,0.4 | 0.571429 | 0.571429 | 0.585403 | 0.585403 | 0.80024 | 0.80024 |
| 0.4,0.5 | 0.836735 | 0.836735 | 0.804898 | 0.804898 | 1.00024 | 1.00024 |
| 0.5,0.5 | 0.918367 | 0.918367 | 0.853662 | 0.853662 | 1.00020 | 1.00020 |

**Table 2.2** Comparison of MTB-DQM solutions for $v$ of Problem 2 with $\mathrm{Re}=80$, $\Delta t=10^{-4}$ with exact solutions and different time $t$

| Mesh | $t=0.1$ | | $t=0.3$ | | $t=0.5$ | |
|---|---|---|---|---|---|---|
| | Num. | Exact | Num. | Exact | Num. | Exact |
| 0.1, 0.1 | -0.02041 | -0.02041 | -0.07321 | -0.07321 | -0.20012 | -0.20012 |
| 0.3, 0.1 | 0.183673 | 0.183673 | 0.170732 | 0.170732 | 0.20004 | 0.20004 |
| 0.2, 0.2 | -0.04082 | -0.04082 | -0.14641 | -0.14641 | -0.40024 | -0.40024 |
| 0.4, 0.2 | 0.163265 | 0.163265 | 0.097527 | 0.097527 | -0.00008 | -0.00008 |
| 0.5, 0.3 | 0.142857 | 0.142857 | 0.024321 | 0.024321 | -0.2002 | -0.2002 |
| 0.1, 0.3 | -0.23980 | -0.2398 | -0.46356 | -0.46356 | -1.00052 | -1.00052 |
| 0.3, 0.4 | -0.18367 | -0.18367 | -0.41479 | -0.41479 | -1.00056 | -1.00056 |
| 0.2, 0.4 | -0.28571 | -0.28571 | -0.53676 | -0.53676 | -1.20064 | -1.20064 |
| 0.4, 0.5 | -0.20408 | -0.20408 | -0.48800 | -0.48800 | -1.20068 | -1.20068 |
| 0.5, 0.5 | -0.10204 | -0.10204 | -0.36603 | -0.36603 | -1.0006 | -1.0006 |

**Table 2.3** $L_2$ , $L_\infty$ errors for $v$ in Problem 2 with $\mathrm{Re}=100$ for various grid sizes at $t=0.01, 0.5,\ \Delta t=0.0001$

| $(N_x, N_y)$ | $t=0.5$ | | $t=0.01$ | |
|---|---|---|---|---|
| | $L_2$ | $L_\infty$ | $L_2$ | $L_\infty$ |
| (4, 4) | 6.23E-08 | 2.40E-08 | 1.48E-09 | 4.83E-10 |
| (8, 8) | 4.21E-09 | 1.21E-09 | 1.11E-10 | 2.27E-11 |
| (17, 17) | 2.66E-10 | 4.36E-11 | 7.55E-12 | 9.03E-13 |
| (32, 32) | 1.61E-11 | 1.40E-12 | 5.10E-13 | 3.13E-14 |
| (44, 44) | 4.13E-12 | 2.64E-13 | 1.60E-13 | 6.71E-15 |
| (64, 64) | 6.28E-13 | 2.28E-14 | 5.43E-14 | 1.80E-15 |

**Table 2.4** $L_2, L_\infty$ errors for $u$ in Problem 2 with $\mathrm{Re}=100$ for various grid sizes and $t=0.01,\ 0.5$ and $\Delta t=0.0001$

| $(N_x, N_y)$ | $t=0.5$ | | $t=0.01$ | |
|---|---|---|---|---|
| | $L_2$ | $L_\infty$ | $L_2$ | $L_\infty$ |
| (4, 4) | 6.23E-08 | 2.40E-08 | 1.48E-09 | 4.83E-10 |
| (8, 8) | 4.21E-09 | 1.21E-09 | 1.11E-10 | 2.27E-11 |
| (17, 17) | 2.66E-10 | 4.36E-11 | 7.55E-12 | 9.03E-13 |
| (32, 32) | 1.61E-11 | 1.40E-12 | 5.10E-13 | 3.13E-14 |
| (44, 44) | 4.13E-12 | 2.64E-13 | 1.60E-13 | 6.71E-15 |
| (64, 64) | 6.28E-13 | 2.28E-14 | 5.43E-14 | 1.80E-15 |

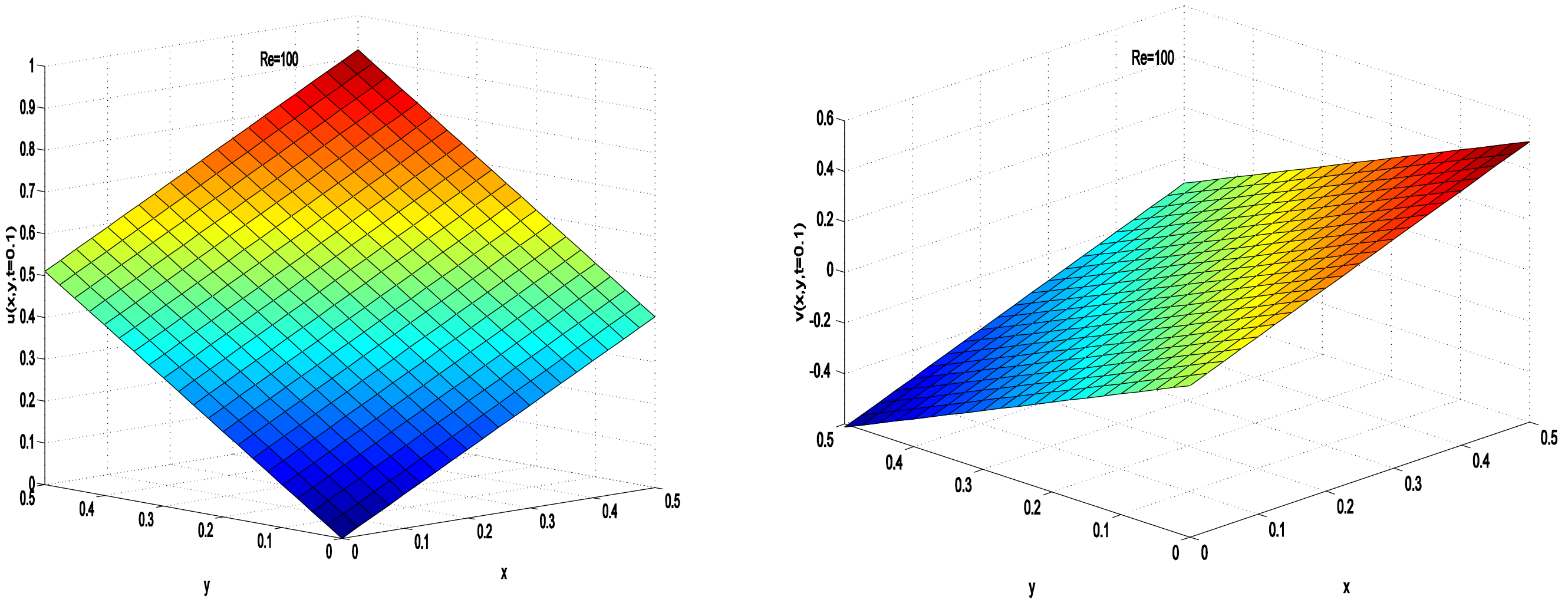

**Figure 5** The physical behavior of MTB-DQM solution of Problem 2 with $\mathrm{Re}=100$ for $u(x, y, t)$ and $v(x, y, t)$ at $t=0.1$

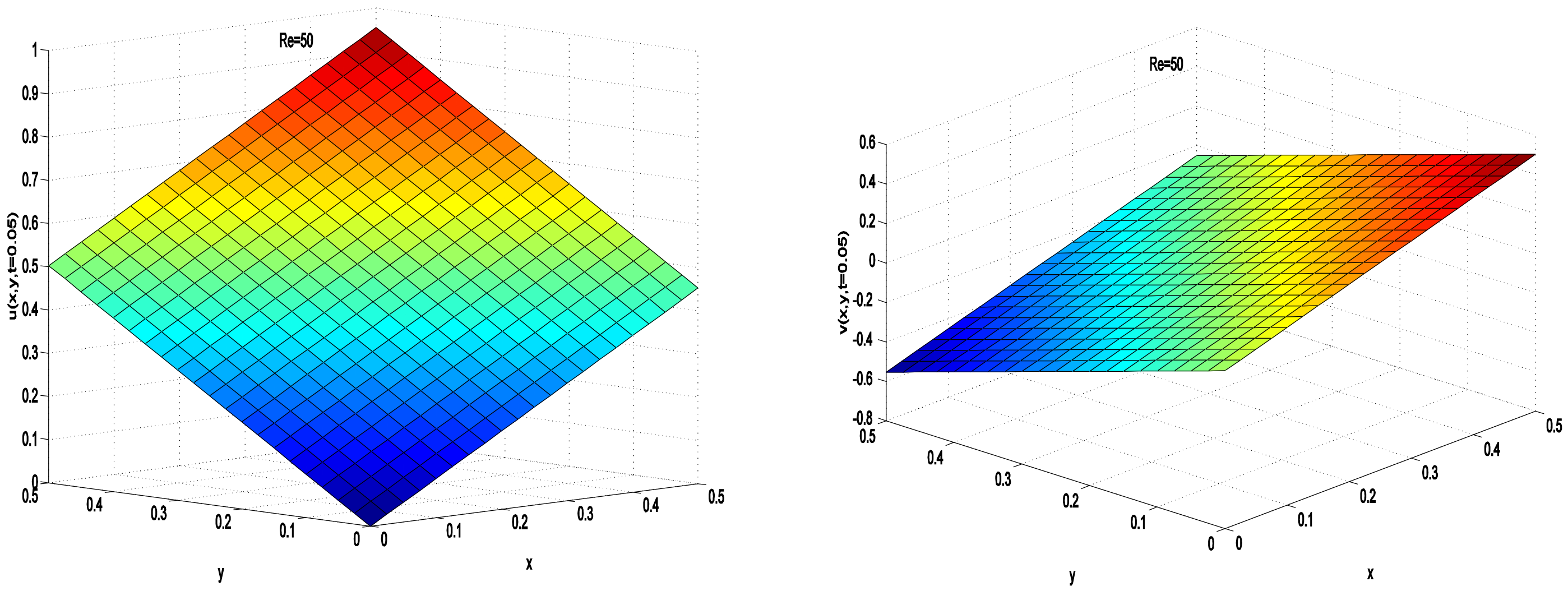

**Figure 6** The physical behavior of MTB-DQM solution of Problem 2 with $\mathrm{Re}=80$ for $u(x,y,t)$ and $v(x,y,t)$ at $t=0.05$

**Problem 3** Consider 2D initial valued coupled viscous Burgers' equations (1.3) in $\Omega=[0,0.5]^2$ with $\phi(x,y)=\sin(\pi x)+\sin(\pi y),\ \psi(x,y)=x+y$ and together with the boundary conditions as follows:

$$\left.\begin{array}{ll} u(0,y,t)=\cos(\pi y), & u(0.5,y,t)=1+\cos(\pi y) \\ v(0,y,t)=y, & v(0.5,y,t)=0.5+y \end{array}\right\},\quad 0\le y\le 0.5,\quad t\ge 0,$$

$$\left.\begin{array}{ll} u(x,0,t)=1+\sin(\pi x), & u(x,0.5,t)=\sin(\pi x) \\ v(x,0,t)=x, & v(x,0.5,t)=x+0.5 \end{array}\right\},\quad 0\le x\le 0.5,\quad t\ge 0.$$

The numerical solutions are obtained at $t=0.1$ for uniform mesh grid $h_x=h_y=0.025$ and $\Delta t=10^{-4}$. The MTB-DQM solutions of $u$ and $v$ components in Problem 3 with $Re=50$ are compared with the solutions in [29, 30, 63] and the exact solutions, for grid size $20\times 20$ and various time levels, see Table 3.1. The physical MTB-DQM solution behavior of $u$ and $v$ components in Problem 3 with $\upsilon=0.01$ at $t=1,2,3$ is depicted in Figure 7.

**Table** 3.1 Comparison of the solutions of Problem 3 with $Re=50$, $(N_x,N_x)=(20,20)$ and $\Delta t=10^{-4}$ at $t=0.625$

| Mesh | $u$ | | | | $v$ | | | |
|---|---|---|---|---|---|---|---|---|
| | MTB-DQM | FDM [30] | Expo-MCB-DQM [63] | Jain-Holla [29] | MTB-DQM | FDM [30] | Expo-MCB-DQM [63] | Jain-Holla[29] |
| 0.1,0.1 | 0.97056 | 0.97146 | 0.97056 | 0.97258 | 0.09842 | 0.09869 | 0.09842 | 0.09773 |
| 0.3,0.1 | 1.15152 | 1.15280 | 1.15152 | 1.16214 | 0.14107 | 0.14158 | 0.14107 | 0.14039 |
| 0.2,0.2 | 0.86244 | 0.86308 | 0.86243 | 0.86281 | 0.16732 | 0.16754 | 0.16732 | 0.16660 |
| 0.4,0.2 | 0.98078 | 0.97985 | 0.98078 | 0.96483 | 0.17223 | 0.17111 | 0.17223 | 0.17397 |
| 0.1,0.3 | 0.66336 | 0.66316 | 0.66335 | 0.66318 | 0.26380 | 0.26378 | 0.26380 | 0.26294 |
| 0.3,0.3 | 0.77226 | 0.77233 | 0.77226 | 0.77030 | 0.22653 | 0.22655 | 0.22653 | 0.22463 |
| 0.2,0.4 | 0.58273 | 0.58181 | 0.58273 | 0.58070 | 0.32935 | 0.32851 | 0.32935 | 0.32402 |
| 0.4,0.4 | 0.76179 | 0.75862 | 0.76179 | 0.74435 | 0.32884 | 0.32502 | 0.32884 | 0.31822 |

**Problem 4** Consider 2D coupled viscous Burgers' equation (1.3) in $\Omega = [0,1]^2$ with exact solution [6]

$$u(x,y,t) = \frac{3}{4} - \frac{1}{4\left(1+\exp\left((-4x+4y-t)\mathrm{Re}/32\right)\right)};\ v(x,y,t) = \frac{3}{4} + \frac{1}{41+\exp\left((-4x+4y-t)\mathrm{Re}/32\right)}$$

where $\phi(x,y)$ and $\psi(x,y)$ on $\Omega$, and $\xi(x,y,t)$, $\zeta(x,y,t)$ on $\partial\Omega$ can be computed from the exact solution for the computational domain $\Omega = [0,1]^2$.

The problem is solved for $\upsilon = 10^{-2}$, $\Delta t = 0.0001$. The computed $L_2$ and $L_\infty$ error norms are compared with the error norm due to the recent schemes mECDQ[45] and Expo-MCB-DQM

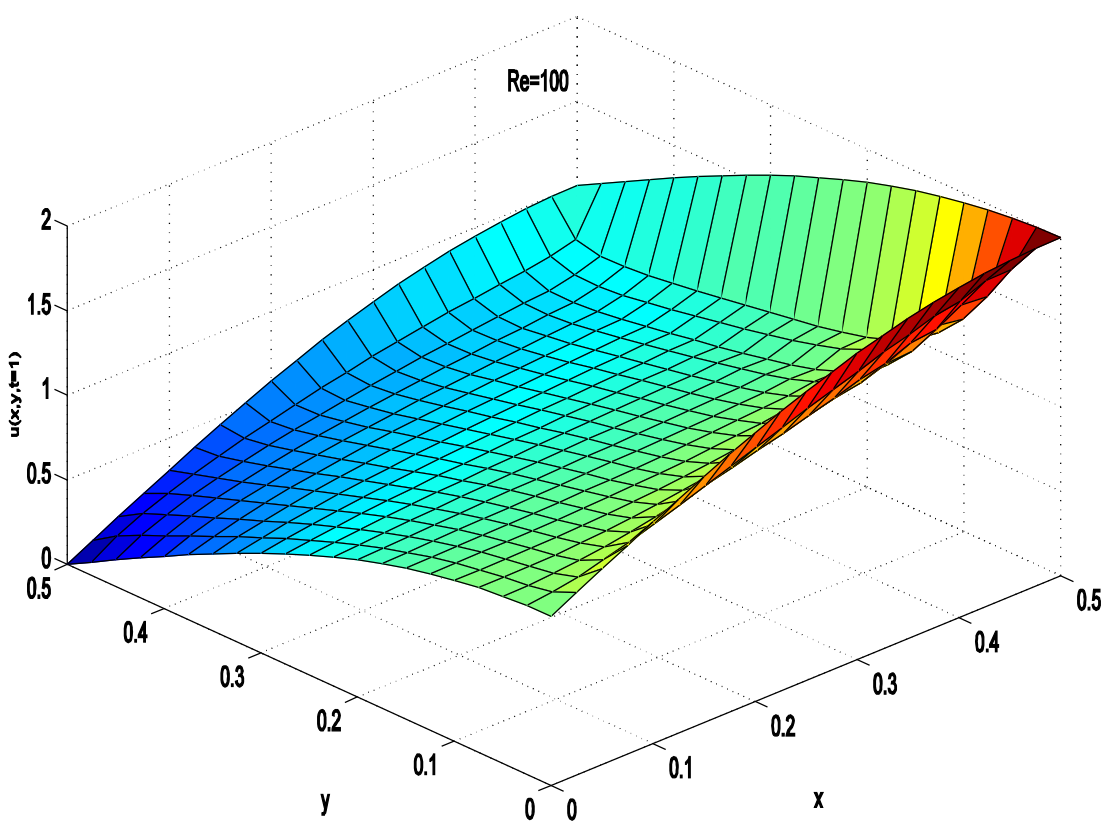

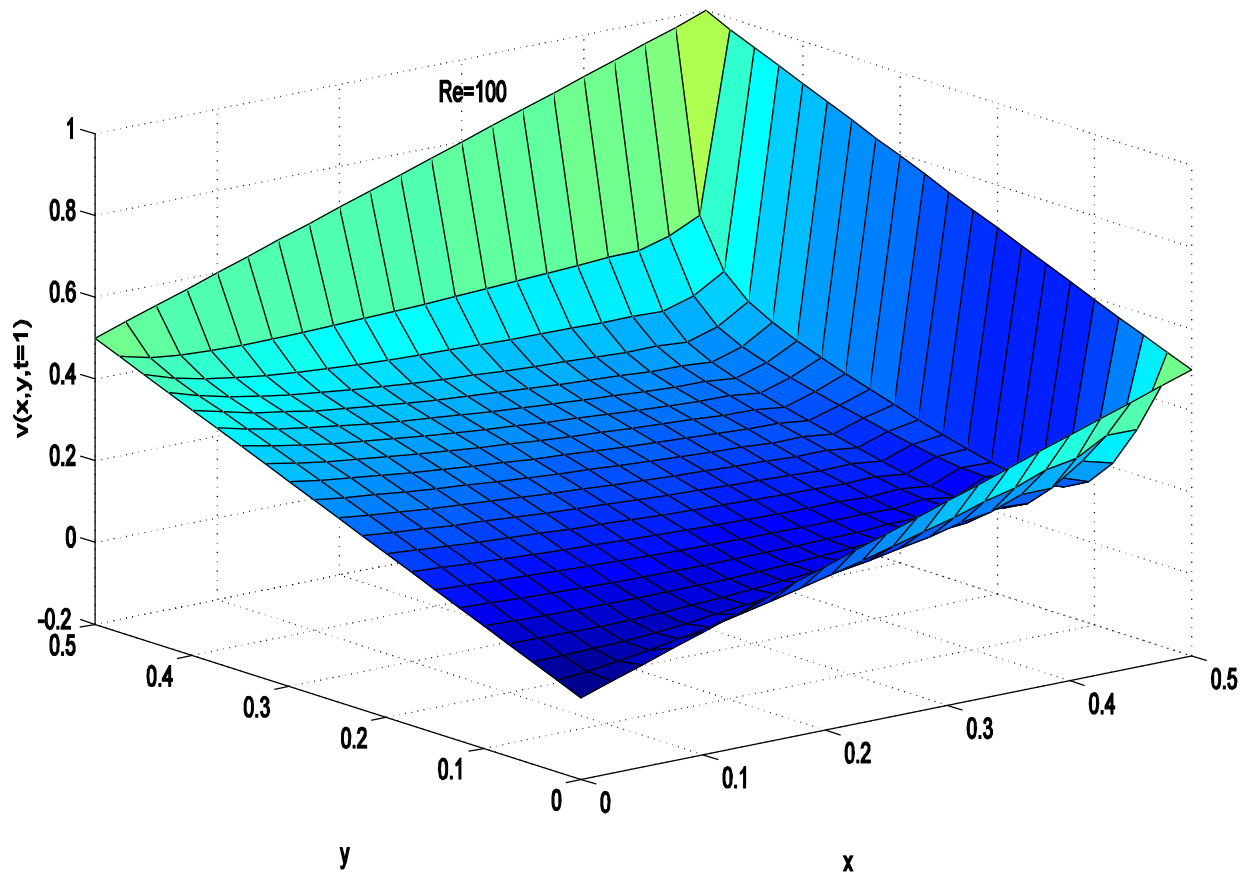

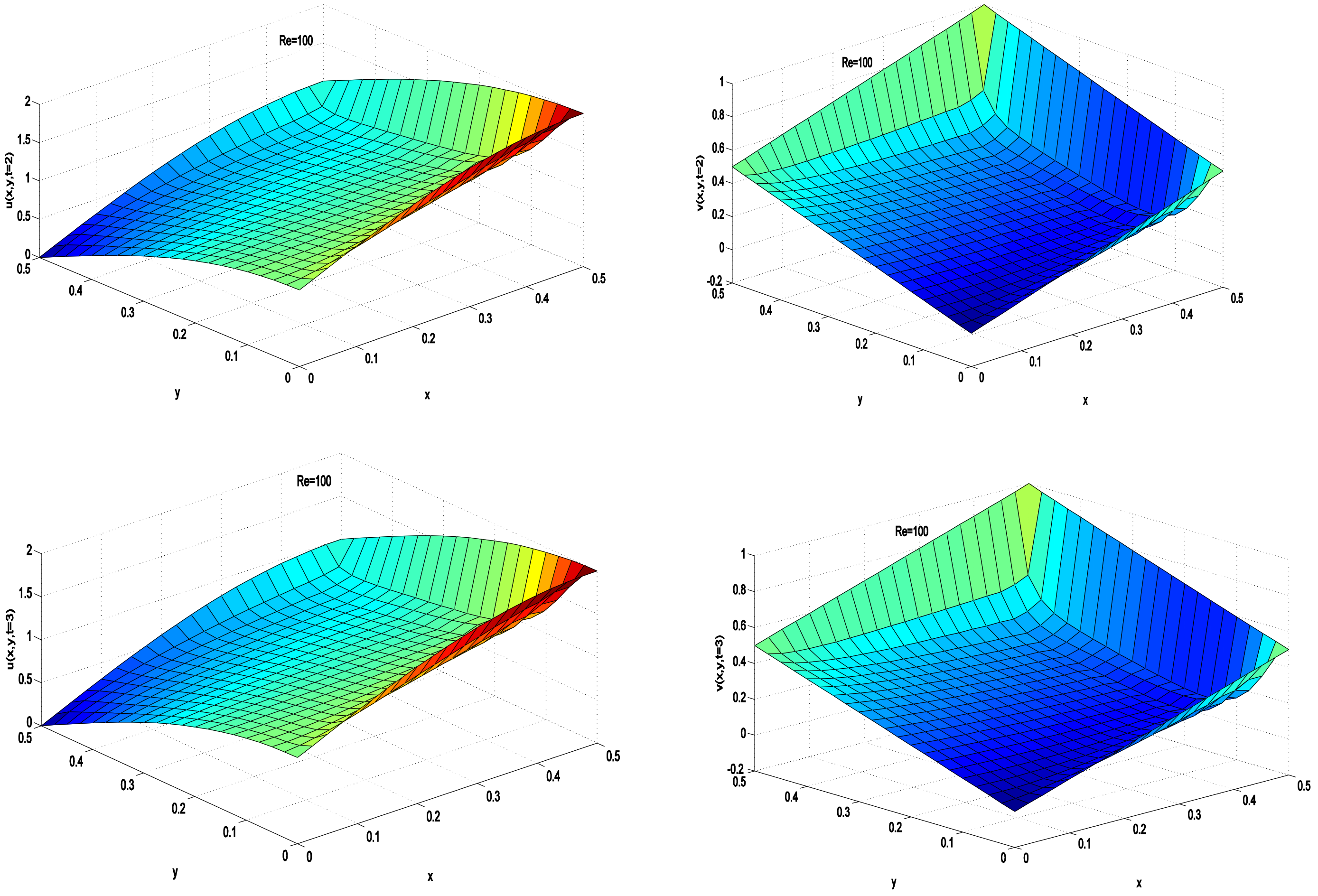

**Figure 7** Physical behavior the numerical solution of Problem 3 for $\mathrm{Re}=100, h_x=h_y=0.025,$ and $\Delta t=10^{-4}$, at $t=3$

[63] and reported in Table 4.1 and Table 4.2. It is evident that the proposed MTB-DQM solutions comparable to [45, 63], and are in good agreement with the exact solutions. The MTB-DQM solutions and exact solutions of $u$ and $v$ at $t=1$ in Problem 4 with $\upsilon=10^{-2}$ and $0.004$ are depicted in Figure 8 and Figure 9, respectively.

**Table 4.1** Comparison of the errors, rate of convergence (R) for $u$ of Problem 4 with $\mathrm{Re}=100, t=0.0001$ at $t=1.0$ and different grid size $(N_x, N_x)$

| $N_x$ | mECDQ [45] $\lambda=0.52$ | | Exp-MCB-DQM [63] $p=10$ | | MTB-DQM | | mECDQ [45] $\lambda=0.52$ | | Exp-MCB- DQM[63] $p=10$ | | MTB-DQM | |
|---|---|---|---|---|---|---|---|---|---|---|---|---|
| | $L_2$ | $R$ | $L_2$ | $R$ | $L_2$ | $R$ | $L_\infty$ | $R$ | $L_\infty$ | $R$ | $L_\infty$ | $R$ |
| *4* | *1.0799E-02* | | *1.5865E-02* | | *1.6451E-02* | | *2.1600E-03* | | *2.3325E-03* | | *2.8958E-03* | |
| *8* | *2.6363E-03* | *2.03* | *1.8037E-03* | *3.13* | *1.9330E-03* | *3.08* | *2.8422E-04* | *2.93* | *1.6816E-04* | *3.79* | *1.9644E-04* | *3.89* |
| *16* | *3.5676E-04* | *2.89* | *3.8329E-04* | *2.23* | *3.9504E-04* | *2.29* | *3.6252E-05* | *3.44* | *1.9610E-05* | *3.10* | *2.0508E-05* | *3.26* |
| *32* | *6.2368E-05* | *2.52* | *8.0461E-05* | *2.52* | *8.1200E-05* | *2.28* | *2.4576E-06* | *3.42* | *2.1967E-06* | *3.15* | *2.2208E-06* | *3.21* |
| *64* | *1.0004E-05* | *2.64* | *1.5355E-05* | *2.38* | *1.5323E-05* | *2.41* | *2.1893E-07* | *3.49* | *2.1795E-07* | *3.33* | *2.1840E-07* | *3.35* |

**Table 4.2** Comparison of the errors, rate of convergence for $v$ of Problem 4 with $\mathrm{Re}=100, t=0.0001$ at $t=1.0$ and different grid size $(N_x, N_x)$

| $N_x$ | mECDQ [45] $\lambda=0.52$ | | Exp-MCB-DQM [63] $p=10$ | | MTB-DQM | | mECDQ [45] $\lambda=0.52$ | | Exp-MCB- DQM[63] $p=10$ | | MTB-DQM | |
|---|---|---|---|---|---|---|---|---|---|---|---|---|
| | $L_2$ | $R$ | $L_2$ | $R$ | $L_2$ | $R$ | $L_\infty$ | $R$ | $L_\infty$ | $R$ | $L_\infty$ | $R$ |
| *4* | *1.0799E-02* | | *1.5865E-02* | | *1.6451E-02* | | *2.1600E-03* | | *2.3325E-03* | | *2.8958E-03* | |
| *8* | *2.6363E-03* | *2.03* | *1.8037E-03* | *3.13* | *1.9330E-03* | *3.08* | *2.8422E-04* | *2.93* | *1.6816E-04* | *3.79* | *1.9644E-04* | *3.89* |
| *16* | *3.5676E-04* | *2.89* | *3.8329E-04* | *2.23* | *3.9504E-04* | *2.29* | *3.6252E-05* | *3.44* | *1.9610E-05* | *3.10* | *2.0508E-05* | *3.26* |
| *32* | *6.2368E-05* | *2.52* | *8.0461E-05* | *2.52* | *8.1200E-05* | *2.28* | *2.4576E-06* | *3.42* | *2.1967E-06* | *3.15* | *2.2208E-06* | *3.20* |
| *64* | *1.0004E-05* | *2.64* | *1.5355E-05* | *2.38* | *1.5323E-05* | *2.41* | *2.1893E-07* | *3.49* | *2.1795E-07* | *3.33* | *2.1840E-07* | *3.35* |

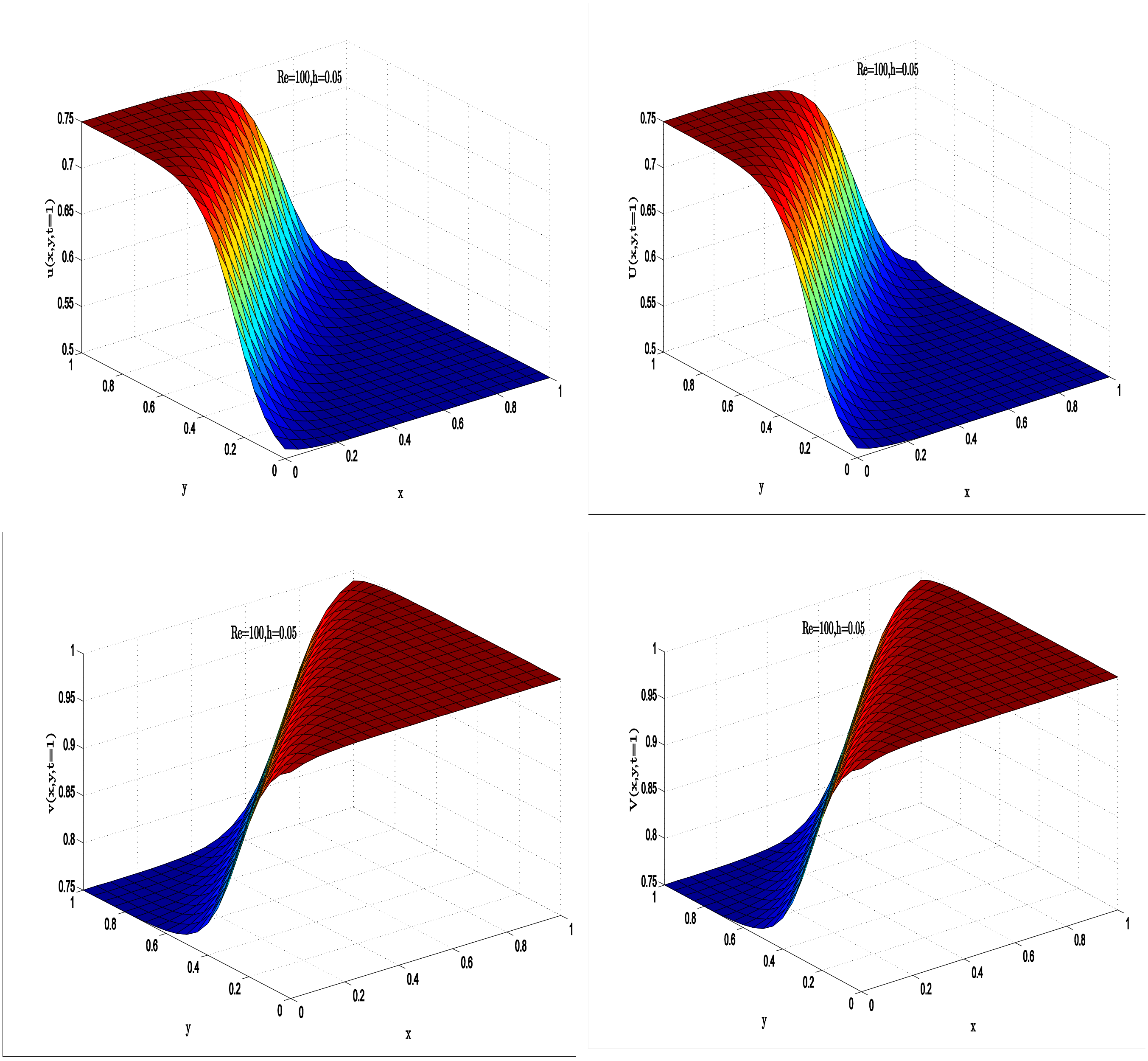

**Figure 8** Comparison of MTB-DQM solution with exact solution of $u$ and $v$ components of Problem 4 with $\upsilon = 0.01, h = 0.05$ and $\Delta t = 0.0001$ at $t = 1$

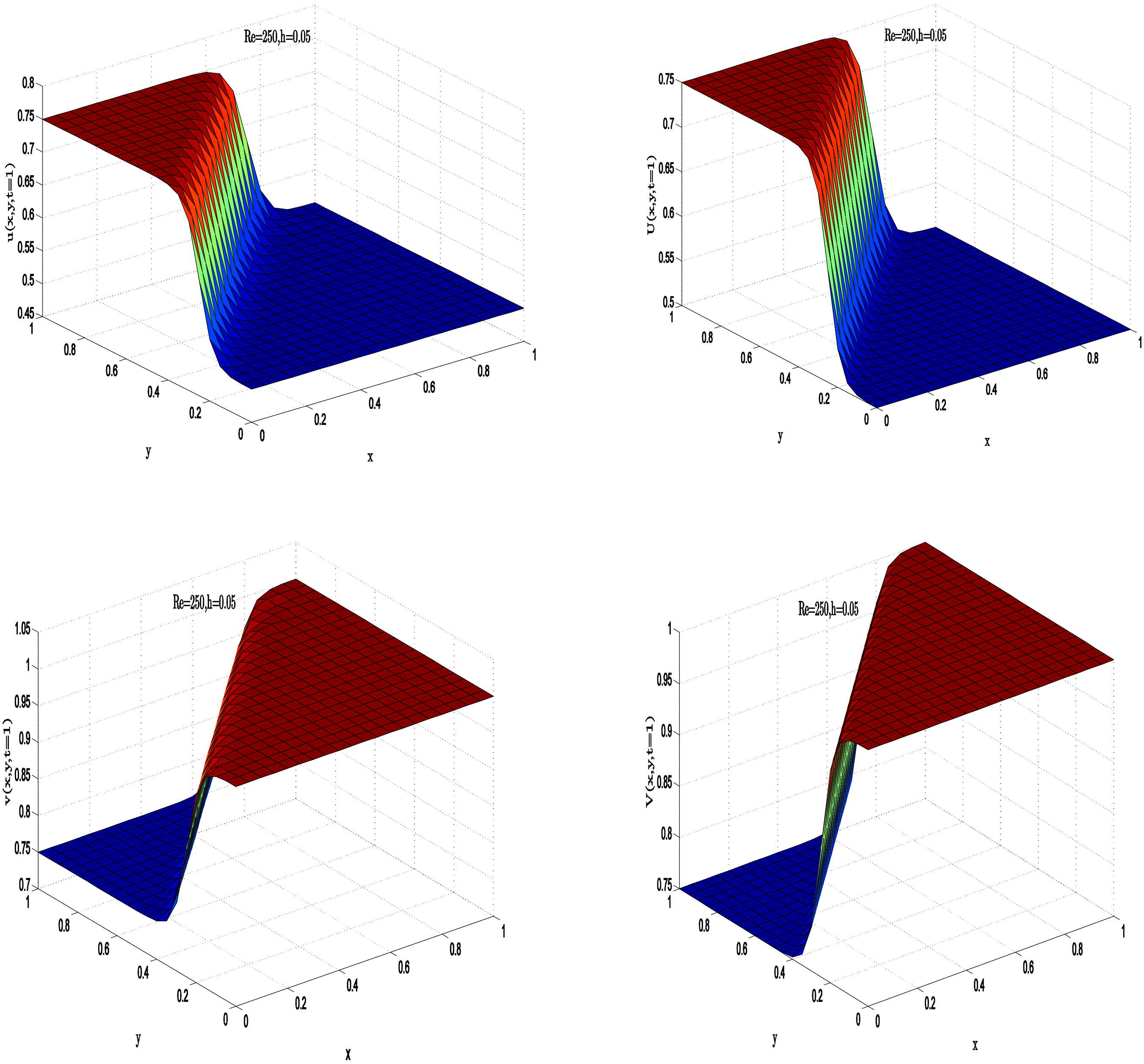

**Figure 9** Comparison of MTB-DQM solution with exact solution of $u$ and $v$ components of Problem 4 with $\mathrm{Re} = 0.004, h = 0.05$ and $\Delta t = 0.0001$ at $t = 1$

## 5. Conclusions

In this paper, a new scheme: "modified trigonometric cubic B-spline differential quadrature method (MTB-DQM)" has been developed for numerical computation of nonlinear partial differential equations. Specially, MTB-DQM has been implemented for the coupled viscous Burgers' equation in one and two dimensions. The MTB-DQM is a new DQM where the weighting coefficients for the spatial derivatives have been computed using modified trigonometric cubic B-splines as set of basis functions. But, we determine the weighting coefficients of the first-order derivative approximations by means of MTB-DQM, and the polynomial based DQM has been adopted to determine the weighting coefficients of the second-order derivative. Thus, MTB-DQM transforms the Burgers 'equation into set of first of order ordinary differential equation

Section 3 confirms that the computational cost of both MCB-DQM [31] and MTB-DQM. In Section 5, the accuracy and efficiency of MTB-DQM has been measured by calculating $L_2$, $L_\infty$ error norms and the rate of convergence, which shows that MTB-DQM generates accurate solutions for Burgers equation in both (1+1) and (2+1) dimensions.

The matrix stability analysis has also been carried out for various grid values which demonstrate that the proposed method is stable for coupled viscous Burgers equation.

## Acknowledgments

P. Kumar is also thankful to Babasaheb Bhimrao Ambedkar University Lucknow India, for financial support to carry out the work.

## References

[1] H. Bateman, Some recentresearches on the motion of fluids, Mon Weather Rev. 43 (1915) 163-170.
[2] J.M. Burgers, Mathematical example illustrating relations occurring in the theory of turbulent fluid motion, Trans. Roy. Neth. Acad. Sci. Amsterdam 17 (1939) 1-53.
[3] J.M. Burgers, A mathematical model illustrating the theory of turbulence, Adv. Appl. Mech. vol. I, Academic Press, New York (1948) 171-199.
[4] S.E Esipov, Coupled Burgers' equations: a model of poly-dispersive sedimentation, Phys. Rev. 52 (1995) 3711-3718.
[5] J. D. Cole, On a quasilinear parabolic equations occurring in aerodynamics, Quart. Appl. Math. 9 (1951) 225-236.
[6] C.A.J. Fletcher, Generating exact solutions of the two dimensional Burgers' equation, Int. J. Numer. Meth. Fluids 3 (1983) 213-216.
[7] J. Nee and J. Duan, Limit set of trajectories of the coupled viscous Burgers' equations, Appl Math Lett 11(1) (1998) 57-61.
[8] E.N. Aksan, Quadratic B-spline finite element method for numerical solution of the Burgers' equation, Appl. Math. Comput. 174 (2006) 884–896.
[9] I.A. Hassanien, A.A. Salama, H.A. Hosham, Fourth-order finite difference method for solving Burgers' equation, Appl. Math. Comput. 170 (2005) 781-800.

[10] M.K. Kadalbajoo, K.K. Sharma, A. Awasthi, A parameter-uniform implicit difference scheme for solving time dependent Burgers' equation, Appl. Math. Comput. 170 (2005) 1365–1393.

[11] W. Liao, An implicit fourth-order compact finite difference scheme for one-dimensional Burgers' equation, Appl. Math. Comput. 206 (2008) 755-764.

[12] H.P. Bhatt, A. Q. M. Khaliq, Fourth-order compact schemes for the numerical simulation of coupled Burgers' equation, Comput. Phy. Commu. 200 (2016) 117 – 138.

[13] R.K. Mohanty, W. Dai, F. Han, Compact operator method of accuracy two in time and four in space for the numerical solution of coupled viscous Burgers' equations, Appl. Math. Comput. 256 (2015) 381–393.

[14] A. Rashid, A.I.Md. Ismail, A Fourier pseudo spectral method for solving coupled viscous Burgers' equations, Comput. Methods Appl. Math. 9 (2009) 412–420.

[15] M.A. Abdou, A.A. Soliman, Variational iteration method for solving Burgers' and coupled Burgers' equations, J. Comput. Appl. Math. 181 (2005) 245–251.

[16] R.C. Mittal, G. Arora, Numerical solution of the coupled viscous Burgers' equation, Commun. Nonlinear Sci. Numer. Simul. 16 (2011) 1304–1313.

[17] I. Dag, D. Irk, A. Sahin, B-Spline collocation methods for numerical solutions of the Burgers' equation, Math. Probl. Eng. 5 (2005) 521-538.

[18] R.C. Mittal, R.K. Jain, Numerical solutions of nonlinear Burgers' equation with modified cubic B-splines collocation method, Appl. Math. Comput. 218 (2012) 7839-7855.

[19] G. Arora, V. Joshi, A computational approach for solution of one dimensional parabolic partial differential equation with application in biological processes, Ain Shams Eng J (2016), http://dx.doi.org/10.1016/j.asej.2016.06.013.

[20] O. Ersoy and I. Dag, The extended B-spline collocation method for numerical solutions of Fisher equation, AIP Conference Proceedings 1648 (2015) 370011. doi: 10.1063/1.4912600.

[21] Shu-Sen Xie, S. Heo, Seokchan Kim, Gyungsoo Woo, Sucheol Yi, Numerical solution of one dimensional Burgers' equation using reproducing kernel function, J. Comput. Appl. Math. 214 (2008) 417-434.

[22] T. Ozis, A. Esen, S. Kutluay, Numerical solution of Burgers' equation by quadratic B-spline finite elements, Appl. Math. Comput. 165 (2005) 237-249.

[23] Huilin Lai, Changfeng Ma, A new lattice Boltzmann model for solving the coupled viscous Burgers equation, Physica A 395 (2014) 445457

[24] . R. C. Mittal and R. Jiwari, Differential quadrature method for numerical solution of coupled viscous Burgers' equations, Int. J. Comput. Methods Engg. Sci. Mechanics, 13 (2) (2012) 88-92.

[25] V.K. Srivastava, M.K. Awasthi, M. Tamsir, A fully implicit finite difference solution to one dimensional coupled Burgers nonlinear equation, Int. J. Math. Sci. 7 (2013) 2328.

[26] Manoj Kumar and SapnaPandit, A composite numerical scheme for the numerical simulation of coupled Burgers equation, Computer Physics Communications 185 (2014) 809-817.

[27] Z.Hongqing,S.Huazhong and D.Meiyu,Numerical solutions of two-dimensional Burgers'equations by discrete Adomian decomposition method,Comp. and Mathematics with App. 60 (2010) 840–848

[28] J. Biazar, H. Aminikhah, Exact and numerical solutions for non-linear Burger's equation by VIM, Math. Comput. Modelling 49 (2009) 1394–1400.

[29] P. C. Jain and D. N. Holla, "Numerical solution of coupled Burgers' equations," Int. J. Numer. Meth. Eng. **12**, 213 (1978).

[30] V. K. Srivastava, S. Singh, and M. K. Awasthi, "Numerical solutions of coupled Burgers' equations by an implicit finite deference scheme," AIP Advances **3**, 082131 (2013).

[31] G. Arora, RC Mittal and B. K. Singh, Numerical Solution of BBM-Burger Equation with Quartic B-spline collocation method, Journal of Engineering Science and Technology, Special Issue 1, 12/2014, 104 - 116.

[32] V. K. Srivastava and B. K. Singh, A Robust finite difference scheme for the numerical solutions of two dimensional time-dependent coupled nonlinear Burgers' equations, International Journal of Applied Mathematics and Mechanics 10(7) (2014) 28-39.

[33] R. Bellman, B. G. Kashef, J. Casti, Differential quadrature: a technique for the rapid solution of nonlinear differential equations, J. Comput. Phy. 10 (1972) 40-52.

[34] J. R. Quan, C.T. Chang, New insights in solving distributed system equations by the quadrature methods-I, Comput. Chem. Eng., 13 (1989) 779-788.

[35] J. R. Quan, C.T. Chang, New insights in solving distributed system equations by the quadrature methods-II, Comput. Chem. Eng., 13 (1989) 1017-1024.

[36] A. Korkmaz, I. Dag, Cubic B-spline differential quadrature methods and stability for Burgers' equation, Eng. Comput. Int. J. Comput. Aided Eng. Software 30 (3) (2013) 320-344.

[37] C. Shu, Y.T. Chew, Fourier expansion-based differential quadrature and its application to Helmholtz eigenvalue problems, Commun. Numer. Methods Eng. 13 (8) (1997) 643-653.

[38] C. Shu, H. Xue, Explicit computation of weighting coefficients in the harmonic differential quadrature J. Sound Vib. 204(3) (1997) 549-555.

[39] A. Korkmaz, I. Dag, Shock wave simulations using sinc differential quadrature method, Eng. Comput. Int. J. Comput. Aided Eng. Software 28 (6) (2011) 654–674.

[40] C. Shu, B. E. Richards, Application of generalized differential quadrature to solve two dimensional incompressible navier-Stokes equations, Int. J. Numer. Meth. Fluids 15 (1992) 791-798.

[41] G. Arora, B. K. Singh, Numerical solution of Burgers' equation with modified cubic B-spline differential quadrature method, Applied Math. Comput. 224 (1) (2013) 166-177.

[42] B. K. Singh and G. Arora, A numerical scheme to solve Fisher-type reaction-diffusion equations, Nonlinear Studies/MESA-MATHEMATICS IN ENGINEERING, SCIENCE AND AEROSPACE 5(2) (2014) 153--164.

[43] B. K. Singh, G. Arora, M. K. Singh, A numerical scheme for the generalized Burgers-Huxley equation, Journal of the Egyptian Mathematical Society (2016) http://dx.doi.org/10.1016/j.joems.2015.11.003.

[44] B. K. Singh and Carlo Bianca, A new numerical approach for the solutions of partial differential equations in three-dimensional space, Appl. Math. Inf. Sci. 10, No. 5, 1-10 (2016).

[45] B.K. Singh, P. Kumar, A novel approach for numerical computation of Burgers' equation (1 + 1) and (2 + 1) dimension, Alexandria Eng. J. (2016) http://dx.doi.org/10.1016/j.aej.2016.08.023

[46] C. Shu, Differential Quadrature and its Application in Engineering, Athenaeum Press Ltd. Great Britain (2000).

[47] W. Lee, Tridiagonal matrices: Thomas algorithm, Scientific Computation, University of Limerick. http://www3.ul.ie/wlee/ms6021_thomas.pdf

[48] S. Gottlieb, D. I. Ketcheson, C. W. Shu, High Order Strong Stability Preserving Time Discretizations, J. Sci. Comput. 38 (2009) 251-289.

[49] J. R. Spiteri, S. J. Ruuth, A new class of optimal high-order strong stability-preserving time-stepping schemes, SIAM J. Numer.Analysis 40 (2) (2002) 469-491.

[50] B. Saka, A. Sahin, I. Dag, B-spline collocation algorithms for numerical solutions of the RLW equation, Numerical Methods for Partial Differential Equations 27(3) (2009) 581-607.

[51] M.K. Jain, Numerical Solution of Differential Equations, 2nd ed., Wiley, New York, NY, 1983.

[52] A. Korkmaz, A.M. Aksoy, I. Dag, Quartic B-spline differential quadrature method, Int. J. Nonlinear Sci. 11 (4) (2011) 403-411.

[53] A. Korkmaz, I. Dag, Polynomial based differential quadrature method for numerical solution of nonlinear Burgers' equation, J. Franklin Inst. 348 (10) (2011) 2863–2875.

[54] A. Bashan, S. B. G. Karakoc, T. Geyikli, Approximation of the KdVB equation by the quintic B-spline differential quadrature method. Kuwait J. Sci. 42 (2) (2015) 67-92.

[55] A. Korkmaz, I. Dag, Quartic and quintic Bspline methods for advection diffusion equation. Applied Mathematics and Computation 274 (2016) 208-219.

[56] A. Korkmaz, and H. K. Akmaz, Numerical Simulations for Transport of Conservative Pollutants. Selcuk Journal of Applied Mathematics 16(1) (2015).

[57] J. Chung, E. Kim, Y. Kim, Asymptotic agreement of moments and higher order contraction in the Burgers equation, J. Diff. Equ. 248 (10) (2010) 2417-2434.

[58] A.H. Salas, Symbolic computation of solutions for a forced Burgers equation, Appl. Math. Comput. 216 (1) (2010) 18–26.

[59] S. Lin, C. Wang, Z. Dai, New exact traveling and non traveling wave solutions for (2+1) dimensional Burgers equation, Appl. Math. Comput. 216 (10) (2010) 3105–3110.

[60] Min Xu, Ren-Hong Wang, Ji-Hong Zhang, Qin Fang, A novel numerical scheme for solving Burgers' equation, Appl. Math. Comput. 217 (2011) 4473-4482.

[61] P. Arminjon, C. Beauchamp, Numerical solution of Burgers' equations in two-space dimensions, Comput. Methods Appl. Mech. Engrg. 19(1979) 351–365.

[62] G.W. Wei, D. S. Zhang, D. J. Kouri, D.K. Hoffman, Distributed approximation functional approach to Burgers' equation in one and two space dimensions, Comput. Phy. Comm. 111 (1998) 93–109.

[63] Mohammad Tamsir, V.K. Srivastava, R. Jiwari, An algorithm based on exponential modified cubic B-spline differential quadrature method for nonlinear Burgers' equation, App. Mathematics and Comp. 290 (2016) 111–124.

[64] Ethan J. Kubatko, Benjamin A. Yeager, David I. Ketcheson, Optimal strong-stability-preserving Runge–Kutta time discretizations for discontinuous Galerkin methods, http://www.davidketcheson.info/assets/papers/dg_ssp_stability.pdf.